\documentclass[12pt,a4paper]{article}
\usepackage{amsmath}
\usepackage{amsfonts}
\usepackage{amssymb}
\usepackage{amsthm}
\usepackage{amscd}
\usepackage{graphicx}
\theoremstyle{plain}
\newtheorem{lemma}{Lemma}[section]

\newtheorem{theorem}{Theorem}[section]

\newtheorem{conjecture}{Conjecture}[section]

\theoremstyle{definition}

\newtheorem{remark}{Remark}[section]

\def\MtDFV{\widetilde{M}}

\def\RtDFV{\widetilde{R}}

\def\RoDFV {R^{\perp}}

\title{Remarks on an inequality involving the normal scalar curvature}

\author{
    \emph{Franki Dillen, Johan Fastenakels and Joeri Van der Veken}
    \thanks{The second author is Research assistant of the Fund for Scientific Research - Flanders (Belgium) (FWO)} \\
    Katholieke Universiteit Leuven\\ Department of Mathematics\\ Celestijnenlaan 200 B - Box 2400\\ BE-3001 Leuven\\
    Belgium\\
    e-mail: \emph{franki.dillen@wis.kuleuven.be} \\
    e-mail: \emph{johan.fastenakels@wis.kuleuven.be}\\
    e-mail: \emph{joeri.vanderveken@wis.kuleuven.be} \\
       }

\date{}

\begin{document}

\maketitle

\begin{abstract}
We study a pointwise inequality for submanifolds in real space forms
involving the scalar curvature, the normal scalar curvature and the
mean curvature. We translate it into an algebraic problem, allowing
us to prove a slightly weaker version of it. We also prove the
conjecture for certain types of submanifolds of ${\mathbb C}^n$.
\end{abstract}

\section{Introduction}

In 1983, Guadelupe and Rodriguez proved the following:

\begin{theorem}[\cite{DFVGR}]
Let $M^2$ be a surface in a real space form $\MtDFV^{2+m}(c)$ of
constant sectional curvature $c$. Denote by $K$ the Gaussian
curvature of $M^2$, by $H$ the mean curvature vector and by
$K^{\perp}$ the normal scalar curvature. Then
$$K \leq \|H\|^2-K^{\perp}+c$$
at every point $p$ of $M^2$, with equality if and only if the
ellipse of curvature at $p$ is a circle.
\end{theorem}
Remark that this is an extension of the well-known inequality
$K\leq\|H\|^2$ for surfaces in ${\mathbb E}^3$.
\medskip

In \cite{DFVDDVV1} the following was conjectured as generalization
of the previous Theorem.

\begin{conjecture}[\cite{DFVDDVV1}]\label{DFVconj}
Let $M^n$ be a submanifold of a real space form $\MtDFV^{n+m}(c)$ of
constant sectional curvature $c$. Denote by $\rho$ the normalized
scalar curvature, by $H$ the mean curvature vector and by
$\rho^{\perp}$ the normalized normal scalar curvature. Then
\begin{equation}\label{DFVineq}
\rho \leq \|H\|^2-\rho^{\perp}+c.
\end{equation}
\end{conjecture}
\medskip

The conjecture was proved for $m=2$ in \cite{DFVDDVV1}, where also
some classification results were obtained in case equality holds in
(\ref{DFVineq}) at every point.

\begin{remark}[Added remark on recent developments]
Nowadays, this conjecture is known as the DDVV-conjecture. Recently
the conjecture was proved for $n=3$ in \cite{DFVCL} and for $m=3$ in
\cite{DFVL1}. In a private communication \cite{DFVL2}, Z. Lu
announced a proof for the general case. Also in the study of
submanifolds attaining equality there is recently substantial
progress : see \cite{DFVDT} and \cite{DFVT}. All these results were
obtained after the finishing of this paper.
\end{remark}

For normally flat submanifolds, in particular for hypersurfaces,
inequality (\ref{DFVineq}) follows from a more general result of
Chen (\cite{DFVC}). In particular, we have for any submanifold $M^n$
of a real space form $\MtDFV^{n+m}(c)$:
\begin{equation}\label{DFV1.1}
\rho \leq\|H\|^2+c.
\end{equation}

For immersions which are invariant with respect to the standard
K\"ahlerian and Sasakian structures on ${\mathbb E}^{2k}$ and
$S^{2k+1}(1)$ the conjecture was proved in \cite{DFVDFV} and for
immersions which are totally real with respect to the nearly
K\"ahler structure on $S^6(1)$ in \cite{DFVDDVV2}.

In section 3 we will translate the conjecture to an algebraic
problem involving symmetric matrices, followed by a proof of a
weaker version. In section 4 we will prove the conjecture for
$H$-umbilical Lagrangian submanifolds of ${\mathbb C}^n\cong{\mathbb
E}^{2n}$, for minimal Lagrangian submanifolds of ${\mathbb
C}^3\cong{\mathbb E}^6$ and for ultra-minimal Lagrangian
submanifolds of ${\mathbb C}^4\cong{\mathbb E}^8$. We remark that
some of these results have been generalized in the meantime by A.
Mihai in \cite{DFVM1}, see \cite{DFVM2} in the present volume. The
reader should be warned however that the notations in \cite{DFVM2}
and in this paper are not always consistent.

\section{Preliminaries}

Let $M^n$ be a Riemannian manifold of dimension $n$ with
Riemann-Christoffel curvature tensor $R$. If $\{e_1,\ldots, e_n\}$
is an orthonormal basis for $T_pM$, then we define the
\emph{normalized scalar curvature} of $M^n$ at $p$ by
\begin{equation}\label{DFVdef scalar curv}
\rho=\frac{2}{n(n-1)}\sum_{i<j=1}^n \langle
R(e_i,e_j)e_j,e_i\rangle.
\end{equation}
\medskip

Now let $\MtDFV^{n+m}$ be another Riemannian manifold with
Riemann-Christoffel curvature tensor $\RtDFV$ and let
$f:M^n\rightarrow\MtDFV^{n+m}$ be an isometric immersion. If $h$ is
the \emph{second fundamental form}, $A_U$  the \emph{shape-operator}
associated to a normal vector field $U$, and  $\RoDFV $ the
curvature tensor of the normal connection, then the equations of
Gauss and Ricci are given by
\begin{equation}\label{DFVGauss eqn}
\langle R(X,Y)Z,T \rangle = \langle \RtDFV(X,Y)Z,T \rangle + \langle
h(X,T),h(Y,Z)\rangle - \langle h(X,Z),h(Y,T) \rangle,
\end{equation}
\begin{equation}\label{DFVRicci eqn}
\langle \RoDFV (X,Y)U,V \rangle = \langle \RtDFV(X,Y)U,V \rangle +
\langle [A_U,A_V]X,Y \rangle,
\end{equation}
for tangent vectors $X$, $Y$, $Z$ and $T$ and normal vectors $U$ and
$V$.
\medskip

Let $\{e_1,\ldots, e_n\}$ be as above and suppose that
$\{u_1,\ldots, u_m\}$ is an orthonormal basis for $T_p^{\perp}M$.
Then we define the \emph{normalized normal scalar curvature} of
$M^n$ at $p$ by
\begin{equation}\label{DFVdef normal scalar curv}
\rho^{\perp}
=\frac{2}{n(n-1)}\sqrt{\sum_{i<j=1}^n\sum_{\alpha<\beta=1}^m\langle
\RoDFV (e_i,e_j)u_{\alpha},u_{\beta}\rangle^2},
\end{equation}
which corresponds to the definition proposed in \cite{DFVDDVV1}.
Another extrinsic curvature invariant that we will use is the
\emph{mean curvature vector} of the submanifold at $p$:
\begin{equation}
H=\frac{1}{n}\mathrm{tr}(h)=\frac{1}{n}\sum_{i=1}^nh(e_i,e_i).
\end{equation}

\section{A translation of the problem}

From now on, we use the following convention: if $A$ and $B$ are
$(n\times n)$-matrices, we define $\langle
A,B\rangle=\mathrm{tr}(A^t\cdot B)$. The associated norm is then
given by $\|A\|^2=\mathrm{tr}(A^t\cdot A)=\sum_{i,j=1}^n
(A_{ij})^2$. The scalar product, and hence the norm are preserved by
orthogonal transformations.

\subsection{The translation} The following theorem reduces the conjecture
to an inequality involving symmetric $(n\times n)$-matrices.
\begin{theorem}\label{DFVtheo 3.1}\label{DFValg}
Conjecture \ref{DFVconj} is true for submanifolds of dimension $n$
and codimension $m$ if for every set $\{B_1,\ldots,B_m\}$ of
symmetric $(n\times n)$-matrices with trace zero the following
inequality holds:
\begin{equation}\label{DFVmatrixineq}
\sum_{\alpha,\beta=1}^m\|[B_{\alpha},B_{\beta}]\|^2\leq\left(\sum_{\alpha=1}^m\|B_{\alpha}\|^2\right)^2.
\end{equation}
\end{theorem}
\emph{Proof.} Let $M^n$ be a submanifold of $\MtDFV^{n+m}(c)$. Take
$p\in M^n$ and suppose that $\{e_1,\ldots,e_n\}$ is an orthonormal
basis for $T_pM$ and that $\{u_1,\ldots,u_m\}$ is an orthonormal
basis for $T_p^{\perp}M$. In summations, Latin indices will always
range from 1 to $n$, whereas Greek indices range from $1$ to $m$.
Further, we use the notations introduced in the previous section.

We define a symmetric $(1,2)$-tensor $b$, taking normal values, by
$$b(X,Y)=h(X,Y)-\langle X,Y\rangle H$$
for all $X,Y\in T_pM$. Remark that
\begin{equation}\label{DFV3.1}
\|b\|^2=\sum_{i,j}\|b(e_i,e_j)\|^2=\sum_{i,j}\|h(e_i,e_j)\|^2-2\sum_{i}\langle
h(e_i,e_i),H\rangle+n\|H\|^2=\|h\|^2-n\|H\|^2.
\end{equation}
Now we define a set $\{B_1,\ldots,B_m\}$ of symmetric operators on
$T_pM$ by
$$\langle B_{\alpha}X,Y\rangle = \langle
b(X,Y),u_{\alpha}\rangle$$ for all $X,Y\in T_pM$. It is clear that
$B_{\alpha}=A_{u_{\alpha}}-\langle H,u_{\alpha}\rangle\mathrm{id}$,
and thus
\begin{equation}\label{DFV3.2}
[B_{\alpha},B_{\beta}]=[A_{u_{\alpha}},A_{u_{\beta}}].
\end{equation}

Using the equation of Gauss (\ref{DFVGauss eqn}) and (\ref{DFV3.1}),
we find
\begin{eqnarray*}
\rho &=& \frac{2}{n(n-1)}\sum_{i<j}\langle R(e_i,e_j)e_j,e_i\rangle\\
     &=& \frac{2}{n(n-1)}\sum_{i<j}\left(c+\langle h(e_i,e_i),h(e_j,e_j)\rangle-\|h(e_i,e_j)\|^2\right)\\
     &=& c+\frac{2}{n(n-1)}\left(\frac{n^2}{2}\|H\|^2-\frac{1}{2}\|h\|^2\right)\\
     &=& c+\frac{n}{n-1}\|H\|^2-\frac{1}{n(n-1)}\left(\|b\|^2+n\|H\|^2\right)\\
     &=& c+\|H\|^2-\frac{1}{n(n-1)}\|b\|^2,
\end{eqnarray*}
and thus
\begin{equation*}
\|H\|^2-\rho
+c=\frac{1}{n(n-1)}\|b\|^2=\frac{1}{n(n-1)}\sum_{\alpha}\|B_{\alpha}\|^2\geq
0.
\end{equation*}
From the equation of Ricci (\ref{DFVRicci eqn}) and (\ref{DFV3.2}),
we get

\begin{equation*}
\rho^{\perp} =
\frac{1}{n(n-1)}\sqrt{\sum_{\alpha,\beta}\|[A_{u_{\alpha}},A_{u_{\beta}}]\|^2}=\frac{1}{n(n-1)}\sqrt{\sum_{\alpha,\beta}\|[B_{\alpha},B_{\beta}]\|^2}.
\end{equation*}

We conclude that
\begin{eqnarray*}
\rho \leq \|H\|^2-\rho^{\perp}+c &\Leftrightarrow& (\rho^{\perp})^2 \leq (\|H\|^2-\rho+c)^2\\
&\Leftrightarrow&
\sum_{\alpha,\beta}\|[B_{\alpha},B_{\beta}]\|^2\leq\left(\sum_{\alpha}\|B_{\alpha}\|^2\right)^2.
\end{eqnarray*}
\hfill $\square$

\begin{remark}
By proving Theorem 2 for $m=2$, we obtain a simple proof of the
conjecture for codimension 2 submanifolds:
$$
\left(\|B_1\|^2 + \|B_2\|^2\right)^2 \geq 4 \|B_1\|^2 \|B_2\|^2 \geq
2\|[B_1, B_2] \|^2, $$ where the second inequality is due to Chern,
do Carmo and Kobayashi \cite{DFVCCK}, see Lemma \ref{DFVlem 3.1}
below.
\end{remark}

\subsection{Proof of a weaker version of the inequality} First, we recall two inequalities.
\begin{lemma}[\cite{DFVCCK}] \label{DFVlem 3.1}
If $B_1$ and $B_2$ are symmetric $(n\times n)$-matrices, then
$$\|[B_1,B_2]\|^2 \leq 2\|B_1\|^2\|B_2\|^2,$$
with equality if and only if $B_1=B_2=0$ or, after a suitable
orthogonal transformation,
\begin{equation}\label{DFVCDK}B_1=\left(\begin{array}{ccccc} 0 & \mu_1 & 0 & \cdots &0\\\mu_1 & 0 & 0 &\cdots &0\\
 0 & 0 & 0 & \cdots & 0\\ \vdots & \vdots & \vdots & \ddots &
\vdots\\ 0&0&0&\cdots&0 \end{array}\right), \quad
B_2=\left(\begin{array}{ccccc} \mu_2 & 0 & 0 & \cdots &0\\0 & -\mu_2 & 0 &\cdots &0\\
 0 & 0 & 0 & \cdots & 0\\ \vdots & \vdots & \vdots & \ddots &
\vdots\\ 0&0&0&\cdots&0 \end{array}\right).\end{equation}
\end{lemma}

\begin{theorem}[\cite{DFVLL}]\label{DFVtheo 3.3}
Let $\{B_1,\ldots,B_m\}$ be a set of symmetric $(n\times
n)$-matrices. Then
\begin{equation*}
\sum_{\alpha,\beta=1}^m\|[B_{\alpha},B_{\beta}]\|^2+\sum_{\alpha,\beta=1}^m
\langle B_{\alpha},B_{\beta}\rangle^2
\leq\frac{3}{2}\left(\sum_{\alpha=1}^n\|B_{\alpha}\|^2\right)^2.
\end{equation*}
\end{theorem}

We will use these inequalities to proof the following, weaker
version of conjecture \ref{DFVconj}:
\begin{theorem}\label{DFVtheo 3.4}
Let $M^n$ be a submanifold of a real space form $\MtDFV^{n+m}$. Then
\begin{itemize}
\item[(i)]$\rho\leq\|H\|^2-\sqrt{\frac{2m-1}{3m-3}}\,\rho^{\perp}+c$,
\item[(ii)]$\rho\leq\|H\|^2-\sqrt{\frac{2}{3}\frac{n^2+n-3}{n^2+n-4}}\,\rho^{\perp}+c$.
\end{itemize}
\end{theorem}
\emph{Proof.} Define the matrices $B_{\alpha}$ as in the proof of
theorem \ref{DFValg}. After a suitable orthogonal transformation, we
may assume that $\langle B_{\alpha},B_{\beta}\rangle=
\|B_{\alpha}\|^2\delta_{\alpha\beta}$. The inequality of
Cauchy-Schwarz yields
$\left(\sum_{\alpha=1}^m\|B_{\alpha}\|^2\right)^2\leq
m\sum_{\alpha=1}^m\|B_{\alpha}\|^4$, and thus
$$\sum_{\alpha\neq\beta=1}^m\|B_{\alpha}\|^2\|B_{\beta}\|^2\leq
(m-1)\sum_{\alpha=1}^m\|B_{\alpha}\|^4.$$ This inequality, together
with lemma \ref{DFVlem 3.1} and theorem \ref{DFVtheo 3.3} gives
\begin{eqnarray}
\frac{3}{2}\left(\sum_{\alpha=1}^n\|B_{\alpha}\|^2\right)^2 &\geq&
\sum_{\alpha,\beta=1}^m\|[B_{\alpha},B_{\beta}]\|^2 +
\sum_{\alpha=1}^n\|B_{\alpha}\|^4\nonumber\\
&\geq& \sum_{\alpha,\beta=1}^m\|[B_{\alpha},B_{\beta}]\|^2 +
\frac{1}{m-1}\left(\sum_{\alpha\neq\beta=1}^m\|B_{\alpha}\|^2\|B_{\beta}\|^2\right)\nonumber\\
&\geq& \sum_{\alpha,\beta=1}^m\|[B_{\alpha},B_{\beta}]\|^2 +
\frac{1}{2(m-1)}\sum_{\alpha,\beta=1}^m\|[B_{\alpha},B_{\beta}]\|^2\nonumber\\
&\geq&
\frac{2m-1}{2m-2}\sum_{\alpha,\beta=1}^m\|[B_{\alpha},B_{\beta}]\|^2.\label{DFV3.4}
\end{eqnarray}
Inequality (\ref{DFV3.4}) implies
$$\frac{3}{2}(\|H\|^2-\rho+c)^2\geq\frac{2m-1}{2m-2}(\rho^{\perp})^2,$$
which yields the first inequality stated in the theorem. To prove
the second one, remark that we may replace $m$ by the dimension of
the image of $b$. The result follows from the observation
$\dim(\mathrm{im}(b))\leq\frac{n(n+1)}{2}-1$.
\hfill$\square$\medskip

\section{Lagrangian submanifolds}
In this section, we prove the conjecture for three families of
Lagrangian submanifolds, namely for $H$-umbilical Lagrangian
submanifolds of ${\mathbb C}^n\cong{\mathbb E}^{2n}$, for minimal
Lagrangian submanifolds of ${\mathbb C}^3\cong{\mathbb E}^6$ and for
ultraminimal Lagrangian submanifolds of ${\mathbb C}^4\cong{\mathbb
E}^8$. Recall that a submanifold $M$ of a K\"ahlerian manifold
$\MtDFV^{2n}$ is called \emph{Lagrangian} if at every point the
almost complex structure $J$ of $\MtDFV^{2n}$ induces an isomorphism
between $T_pM$ and $T_p^{\perp}M$. In particular $\dim(M)=n$. The
second fundamental form satisfies the following symmetry property:
\begin{equation}\label{DFV5.00}
\langle h(X,Y),JZ\rangle = \langle h(X,Z),JY\rangle,
\end{equation}
for $X,Y,Z\in T_pM$.

\subsection{$H$-umbilical Lagrangian immersions in ${\mathbb C}^n$}
It was proven in \cite{DFVCO} that there are no totally umbilical
Lagrangian submanifolds in complex space forms, except totally
geodesic ones. $H$-umbilical Lagrangian submanifolds are introduced
in \cite{DFVC2} as the `simplest' Lagrangian submanifolds next to
totally geodesic ones. Their second fundamental form satisfies
\begin{equation}\label{DFV5.0}
\begin{array}{ll}
h(E_1,E_1)=\lambda JE_1, & h(E_2,E_2)=\ldots=h(E_n,E_n)=\mu JE_1,\\
h(E_1,E_j)=\mu JE_j, & h(E_j,E_k)=0\mbox{ for }j,k\in\{2,\ldots,n\},
j\neq k.
\end{array}
\end{equation}
for some suitable functions $\lambda$ and $\mu$ and a suitable
orthonormal local frame field $\{E_1,\ldots,E_n\}$ on $M^n$.
\medskip

We prove the following:
\begin{theorem} Let $M^n$ be a $H$-umbilical Lagrangian immersion in
${\mathbb C}^n\cong{\mathbb E}^{2n}$. Then
$$\rho \leq \|H\|^2-\rho^{\perp},$$
with equality at every point if and only if $M^n$ is totally
geodesic.
\end{theorem}
\emph{Proof.} From (\ref{DFV5.0}) the form of the shape-operators is
easily deduced. We now use theorem \ref{DFVtheo 3.1}. Defining the
matrices $B_{\alpha}$ as in the proof of that theorem, we easily see
that
\begin{eqnarray*}
\sum_{\alpha,\beta=1}^n\|[B_{\alpha},B_{\beta}]\|^2 &=& 2(n-1)\mu^2\left((n-2)\mu^2+2(\lambda-\mu)^2\right),\\
\left(\sum_{\alpha=1}^n \|B_{\alpha}\|^2\right)^2 &=&
(n-1)^2\left(\frac{1}{n}(\lambda-\mu)^2+2\mu^2\right)^2,
\end{eqnarray*}
such that
\begin{equation*}
\sum_{\alpha,\beta=1}^n\|[B_{\alpha},B_{\beta}]\|\leq\left(\sum_{\alpha=1}^n\|B_{\alpha}\|^2\right)^2
\Leftrightarrow
2n\mu^4-\frac{4}{n}\mu^2(\lambda-\mu)^2+\frac{n-1}{n^2}(\lambda-\mu)^4\geq
0 .
\end{equation*}
The last inequality is satisfied for every $\lambda$ and $\mu$ since
the bilinear form $2nx^2-\frac{4}{n}xy+\frac{n-1}{n^2}y^2$ is
positive definite. \hfill$\square$

\subsection{Minimal Lagrangian submanifolds of ${\mathbb C}^3$}
\begin{theorem}\label{DFVtheo 5.2}
Let $M^3$ be a minimal Lagrangian submanifold of ${\mathbb C}^3$.
Then
$$\rho\leq-\rho^{\perp}$$
and equality holds at a point $p$ if and only if there exists an
orthonormal basis $\{e_1,e_2,e_3\}$ of $T_pM$ such that
\begin{equation}\label{DFVS3}
A_{Je_1}=\left(\begin{array}{ccc}a&0&0\\0&-a&0\\0&0&0\end{array}\right),
\quad
A_{Je_2}=\left(\begin{array}{ccc}0&-a&0\\-a&0&0\\0&0&0\end{array}\right),
\quad A_{Je_3}=0,
\end{equation}
with respect to this basis. If equality holds at every point of a
minimal Lagrangian submanifold of ${\mathbb C}^3$, then $M^3$ is
either a cylinder on complex curve in ${\mathbb C}^2$ (with respect
to a different complex structure) or a ``twisted special Lagrangian
cone", both in the sense of \cite{DFVB}.
\end{theorem}

\emph{Proof.} Let $M^3$ be a minimal Lagrangian submanifold of
${\mathbb C}^3$. Take $p\in M^3$ and consider the function
$$f:\left\{X\in T_pM\ |\ \|X\|=1\right\}\rightarrow{\mathbb R}:X\mapsto\langle h(X,X),JX\rangle.$$
Take $e_1\in T_pM$ such that $f$ attains its maximum value in $e_1$.
Then $\langle h(e_1,e_1),JY\rangle=0$ for every $Y\perp e_1$. Using
(\ref{DFV5.00}), this implies that $e_1$ is an eigenvector of
$A_{Je_1}$. Choosing $e_2$ and $e_3$ such that $\{e_1,e_2,e_3\}$ is
an orthonormal basis for $T_pM$ which diagonalizes $A_{Je_1}$, we
have that the shape-operators take the following form:
\begin{equation*}
A_{Je_1}=\left(\begin{array}{ccc}a+b&0&0\\0&-a&0\\0&0&-b\end{array}\right),
\quad
A_{Je_2}=\left(\begin{array}{ccc}0&-a&0\\-a&c&-d\\0&-d&-c\end{array}\right),
\quad
A_{Je_3}=\left(\begin{array}{ccc}0&0&-b\\0&-d&-c\\-b&-c&d\end{array}\right).
\end{equation*}

We now compute $\rho$, using Gauss' equation:
\begin{eqnarray*}
3\rho &=& \sum_{i<j=1}^3\langle R(e_i,e_j)e_j,e_i\rangle\\
&=& \sum_{i<j=1}^3\left(\langle h(e_i,e_i),h(e_j,e_j)\rangle-\langle h(e_i,e_j),h(e_i,e_j)\rangle\right)\\
&=& (-2a^2-ab)+(-2b^2-ab)+(ab-2c^2-2d^2)\\
&=& -2(a^2+b^2+c^2+d^2)-ab.
\end{eqnarray*}
The computation of $\rho^{\perp}$ using Ricci's equation is
completely analoguous to that in \cite{DFVDDVV2}, yielding
\begin{eqnarray*}
9(\rho^{\perp})^2 &=& \sum_{\alpha<\beta=1}^3\sum_{i<j=1}^3\langle \RoDFV (e_i,e_j)Je_{\alpha},Je_{\beta}\rangle^2\\
&=&  \frac{1}{2}\sum_{\alpha<\beta=1}^3\|[A_{Je_{\alpha}},A_{Je_{\beta}}]\|^2\\
&=& 4(a^4+b^4+c^4+d^4)+4a^3b+4ab^3\\
& & +3a^2b^2+2a^2c^2+2a^2d^2+2b^2c^2+2b^2d^2+8c^2d^2-8abc^2-8abd^2.
\end{eqnarray*}

Using the same argument as in \cite{DFVDDVV2}, we obtain that
$9(\rho^{\perp})^2\leq(3\rho)^2$, which implies the inequality
stated in the theorem, since $\rho\leq0$ from (\ref{DFV1.1}).
Equality holds if and only if $c=d=0$ and $ab=0$. By, if necessary,
changing the role of $e_2$ and $e_3$, we obtain the result.

For proving the statement on the equality case, it suffices to
remark that when the shape operator has the form (\ref{DFVS3}), then
the cubic form $\langle h(X,Y), JZ\rangle$ has $S_3$-symmetry in the
sense of \cite{DFVB} and therefore the classification following from
the classification in \cite{DFVB}.

\hfill$\square$
\medskip

We can extend the previous theorem to 3-dimensional Lagrangian
submanifolds of complex space forms. For a complex space form of
constant holomorphic sectional curvature $4c$, the curvature tensor
takes the form
\begin{equation*}
\RtDFV(X,Y)Z=c\left(\langle Y,Z\rangle X-\langle X,Z\rangle
Y+\langle JY,Z\rangle JX-\langle JX,Z\rangle JY-2\langle JX,Y\rangle
JZ\right).
\end{equation*}
This implies that for a Lagrangian immersion in such a space, the
equations of Gauss and Ricci read respectively:
\begin{equation*}
\langle R(X,Y)Z,T \rangle = c\left(\langle Y,Z\rangle\langle
X,T\rangle-\langle X,Z\rangle\langle Y,T\rangle\right) + \langle
h(X,T),h(Y,Z)\rangle - \langle h(X,Z),h(Y,T) \rangle,
\end{equation*}
\begin{equation*}
\langle \RoDFV (X,Y)U,V \rangle = c\left(\langle JY,U\rangle\langle
JX,V\rangle-\langle JX,U\rangle\langle JY,V\rangle\right) + \langle
[A_U,A_V]X,Y \rangle.
\end{equation*}
An analoguous computation as in the proof of the previous theorem
now yields the following:
\begin{theorem}
Let $M^3$ be a minimal Lagrangian submanifold of a complex space
form of constant holomorphic sectional curvature $4c$. Then
$$(\rho^{\perp})^2\leq(\rho-c)^2-2c(\rho-c)+\frac{c^2}{3}$$
and equality holds at a point $p$ if and only if there exists an
orthonormal basis $\{e_1,e_2,e_3\}$ of $T_pM$ such that
\begin{equation}
A_{Je_1}=\left(\begin{array}{ccc}a&0&0\\0&-a&0\\0&0&0\end{array}\right),
\quad
A_{Je_2}=\left(\begin{array}{ccc}0&-a&0\\-a&0&0\\0&0&0\end{array}\right),
\quad A_{Je_3}=0,
\end{equation}
with respect to this basis.
\end{theorem}
In \cite{DFVDFV} an analoguous inequality relating $\rho$ and
$\rho^{\perp}$ is obtained for complex submanifolds of complex space
forms.

\subsection{Ultra-minimal Lagrangian submanifolds of ${\mathbb C}^4$} A
submanifold $M^n$ of a Riemannian manifold $\MtDFV^{n+m}$ is called
ultra-minimal if around each point $p\in M^n$ there exist a local
orthonormal tangent frame and a local orthonormal normal frame, such
that the shape operators take the form
$$A_{U_\alpha}=\left(\begin{array}{cccccc} A_1^{\alpha} &\cdots & 0 & 0 & \cdots & 0\\
\vdots & \ddots & \vdots & \vdots  & \ddots & \vdots\\
0 & \cdots & A_k^{\alpha} & 0 & \cdots & 0\\
0 & \cdots & 0 & 0 & \cdots & 0\\
\vdots & \ddots & \vdots & \vdots & \ddots & \vdots\\
0 & \cdots & 0 & 0 & \cdots & 0
\end{array}\right),$$
where $A_j^{\alpha}$ is an symmetric $(n_j\times n_j)$-matrix, with
$\mathrm{tr}(A_j^{\alpha})=0$, and $n_1\neq n$.

\begin{theorem}
Let $M^4$ be an ultra-minimal submanifold of ${\mathbb C}^4$. Then
$$\rho\leq-\rho^{\perp},$$
and equality holds at a point $p$ if and only if there exists an
orthonormal basis $\{e_1,e_2,e_3,e_4\}$ of $T_pM$ such that
\begin{equation}\label{DFV5.1}
A_{Je_1}=\left(\begin{array}{cccc}a&b&0&0\\b&-a&0&0\\0&0&0&0\\0&0&0&0\end{array}\right),
\quad
A_{Je_2}=\left(\begin{array}{cccc}b&-a&0&0\\-a&-b&0&0\\0&0&0&0\\0&0&0&0\end{array}\right),
\quad A_{Je_3}=0, \quad A_{Je_4}=0,
\end{equation}
with respect to this basis.
\end{theorem}

\emph{Proof.} Since $M^4$ is ultra-minimal, there are two cases to
consider, namely $n_1=n_2=2$ and $n_1=3$, $n_2=1$.

In the first case, using the symmety conditions for Lagrangian
immersions, we obtain that
\begin{eqnarray*}
A_{Je_1}=\left(\begin{array}{cccc}a&b&0&0\\b&-a&0&0\\0&0&0&0\\0&0&0&0\end{array}\right),
& &
A_{Je_2}=\left(\begin{array}{cccc}b&-a&0&0\\-a&-b&0&0\\0&0&0&0\\0&0&0&0\end{array}\right),\\
A_{Je_3}=\left(\begin{array}{cccc}0&0&0&0\\0&0&0&0\\0&0&c&d\\0&0&d&-c\end{array}\right),
& &
A_{Je_4}=\left(\begin{array}{cccc}0&0&0&0\\0&0&0&0\\0&0&c&-d\\0&0&-d&-c\end{array}\right).
\end{eqnarray*}
Using Ricci's equation, one can verify that
$$36(\rho^{\perp})^2=\frac{1}{2}\sum_{\alpha<\beta=1}^4\|[A_{Je_{\alpha}},A_{Je_{\beta}}]\|^2=4\left((a^2+b^2)^2+(c^2+d^2)^2\right)$$
and from Gauss' equation
$$6\rho=-2(a^2+b^2+c^2+d^2),$$ thus we have $\rho\leq-\rho^{\perp}$,
with equality if and only if $a=b=0$ or $c=d=0$.

In the second case, the ultra-minimality condition yields that
$A_2^{\alpha}=0$ for $\alpha=1,2,3,4$ and hence the problem reduces
to the one solved in Theorem \ref{DFVtheo 5.2}. We obtain
$\rho\leq-\rho^{\perp}$, with equality if and only if the shape
operators take the form (\ref{DFV5.1}), with $b=0$.
\hfill$\square$\medskip

%

\end{document}